\documentclass{amsart}
\usepackage{amssymb} 
\usepackage{amsmath} 
\usepackage{amscd}
\usepackage[matrix,arrow]{xy}

\DeclareMathOperator{\Res}{Res}
\DeclareMathOperator{\rank}{rank}
\DeclareMathOperator{\Pic}{Pic}
\DeclareMathOperator{\red}{red}
\newcommand{\vv}{{\upsilon}}

\DeclareMathOperator{\ord}{ord}

\newcommand{\crt}{\sqrt[3]{2}}
\newcommand{\Q}{{\mathbb Q}}
\newcommand{\Z}{{\mathbb Z}}
\newcommand{\F}{{\mathbb F}}

\newcommand{\cB}{\mathcal{B}}
\newcommand{\cK}{\mathcal{K}}

\newcommand{\cQ}{\mathcal{Q}}

\newcommand{\cC}{\mathcal C}

\newcommand{\cW}{\mathcal W}
\newcommand{\OO}{{\mathcal O}}

\begin {document}

\newtheorem{thm}{Theorem}
\newtheorem{lem}{Lemma}[section]

\theoremstyle{definition}

\theoremstyle{remark}

\title[Explicit Chabauty]{Explicit Chabauty over Number Fields}
\author{Samir Siksek}
\address{Mathematics Institute\\
	University of Warwick\\
	Coventry\\
	CV4 7AL \\
	United Kingdom}

\email{s.siksek@warwick.ac.uk}
\date{\today}
\thanks{The author is supported by an EPSRC Leadership Fellowship}


\begin {abstract}
Let $C$ be a smooth projective absolutely irreducible curve
of genus $g \geq 2$ over a number field $K$ of degree $d$,
and denote its Jacobian by $J$.
Denote the Mordell--Weil rank of $J(K)$ by $r$. 
We give an explicit and practical Chabauty-style criterion for showing that
a given subset $\cK \subseteq C(K)$ is in fact equal to
$C(K)$. 
This criterion is likely to be successful if $r \leq d(g-1)$.
We also show that the only solutions to the equation $x^2+y^3=z^{10}$
in coprime non-zero integers is $(x,y,z)=(\pm 3, -2, \pm 1)$. This  is
achieved by reducing the problem to the determination of $K$-rational points
on several genus $2$ curves where $K=\Q$ or $\Q(\sqrt[3]{2})$, and applying the
method of this paper. 
\end {abstract}
\maketitle

\section{Introduction}
Let $C$ be a smooth projective absolutely irreducible 
curve of genus $g \geq 2$ defined
over a number field $K$, and write $J$ for the Jacobian of $C$. 
Suppose that the rank of the Mordell--Weil group $J(K)$
is at most $g-1$. In a pioneering paper, Chabauty \cite{Chabauty}
proved the finiteness of the set of $K$-rational points on $C$.
This has since been superceded by Faltings' proof of the Mordell
conjecture \cite{Faltings1} which gives the finiteness of $C(K)$ without any
assumption on the rank of $J(K)$. Chabauty's approach,
where applicable, does however have two considerable advantages:

\smallskip

\noindent {\bf (a)} The first is
that Chabauty can be refined to give explicit 
bounds for the cardinality 
of $C(K)$, as shown by Coleman \cite{Co}.
Coleman's bounds
are realistic, and occasionally even sharp; see for example
\cite{Grant}, \cite{FlynnCo}. Coleman's approach has
been adapted to give bounds (assuming some reasonable conditions) 
for the number of solutions
of Thue equations \cite{LT}, the number of rational points on
Fermat curves \cite{Mc1}, \cite{Mc2}, the number of points on
curves of the form $y^2=x^5+A$ \cite{Stoll06},
and the number of rational points on twists of a given curve \cite{Stollch}.

\smallskip

\noindent {\bf (b)} The second is that the Chabauty--Coleman strategy can
often be adapted to compute $C(K)$,
as in \cite{Br1}, \cite{Br2}, \cite{Flynn}, \cite{FW1}, \cite{FW2}, \cite{MP}, 
 \cite{We}, and even the $K$-rational points on the symmetric powers
of $C$ \cite{Sik09}.

\smallskip

This paper is inspired by a 
talk~\footnote{http://msri.org/publications/ln/msri/2000/arithgeo/wetherell/1/banner/01.html} 
given by Joseph Wetherell at MSRI on December 11, 2000. In that talk Wetherell suggested
that it should be possible to adapt the Chabauty strategy to compute the 
set of $K$-rational points on $C$ provided the rank $r$ 
of the Mordell--Weil group
$J(K)$ satisfies $r \leq d(g-1)$, where $d=[K:\Q]$. 
Wetherell has never published
details of his method which we believe is similar to the one we give here. 

In this paper we give a practical Chabauty-style method for determining
$C(K)$ which should succeed if the inequality $r \leq d(g-1)$
holds (but see the discussion at the end of 
Section~\ref{sec:heu}).
We suppose that we have been supplied with a basis
$D_1,\dots,D_r$
for a subgroup of $J(K)$ of full-rank and hence finite
index---the elements of this basis 
are represented as degree $0$ divisors on $C$ (modulo
linear equivalence). Obtaining a basis for a subgroup
of full-rank is often the happy outcome of a successful descent
calculation (see for example 
\cite{CF}, \cite{Flynnd}, 
\cite{PS}, \cite{Sch},  \cite{SchW}, 
\cite{StollI}, \cite{Stoll}, \cite{StollII}). 
Obtaining
a basis for the full Mordell--Weil group is often time consuming for
genus $2$ curves (\cite{Flynnh}, \cite{FS}, \cite{Stollh1}, \cite{Stollh2})
 and simply not feasible in the present state
of knowledge for curves of genus $\geq 3$. We also assume the knowledge
of at least one rational point $P_0 \in C(K)$. If a search for rational
points on $C$ does not reveal any points, then experience suggests
that $C(K)=\emptyset$, and the methods of \cite{BS}, \cite{BS2}
are likely to prove this.

This paper is organized as follows. Section~\ref{sec:heu} gives
a heuristic explanation of why Chabauty's approach should be
applicable when the rank $r$ of the Mordell--Weil group
satisfies $r \leq g(d-1)$. Section~\ref{sec:prelim} gives a quick
summary of basic facts regarding $\vv$-adic integration
on curves and Jacobians. In Section~\ref{sec:chab},
for  $Q \in C(K)$ and a rational prime $p$,
we define a certain neighbourhood of $Q$ in $\prod_{\vv \mid p} C(K_\vv)$
that we call the $p$-unit ball around $Q$,
and give a Chabauty-style criterion for $Q$ 
to be the unique $K$-rational point belonging
to this neighbourhood. In Section~\ref{sec:MW} we explain how
to combine our Chabauty criterion with the Mordell--Weil sieve,
and deduce a practical criterion for a given set $\cK \subseteq C(K)$
to be equal to $C(K)$.
 In Section~\ref{sec:Fermat} we 
use our method to prove the following theorem. 
\begin{thm}\label{thm:Ferm}
The only solutions to the equation
\begin{equation}\label{eqn:2310}
x^2+y^3=z^{10}
\end{equation}
in coprime integers $x$, $y$, $z$ are
\[
(\pm 3, -2, \pm 1), \quad (\pm 1,0 , \pm 1), 
\quad (\pm 1, -1, 0), 
\quad (0,1, \pm 1).
\]
\end{thm}
We note that Dahmen \cite[Chapter 3.3.2]{Dahmen} has solved the 
equation $x^2+z^{10}=y^3$ 
using Galois representations and level lowering. 
We have been unable to solve \eqref{eqn:2310} by using 
Galois representations;
the difficulty arises from the additional \lq non-trivial\rq\ solution
$(x,y,z)=(\pm 3, -2, \pm 1)$ which is not present for the equation $x^2+z^{10}=y^3$.
We solve \eqref{eqn:2310} by reducing the problem to determining the 
$K$-rational points
on several genus $2$ curves where $K$ is either $\Q$ or $\Q(\crt)$.
For all these genus $2$ curves the inequality $r \leq d(g-1)$ is
satisfied and we are able to determine the $K$-rational points
using the method of this paper.

Recently, David Brown has given \cite{DB} an independent and entirely 
different proof of Theorem~\ref{thm:Ferm}. Brown's method is rather
intricate, and makes use of elliptic curve Chabauty, mod $5$ level lowering
and number field enumeration.

\medskip

We would like to thank Tim Dokchitser for useful discussions
and Sander Dahmen for corrections.

\section{A Heuristic Explanation of Wetherell's Idea}\label{sec:heu}

In this section we explain the heuristic idea behind Chabauty's method
and then how the heuristic can be modified for curves over number fields.
Let $C$ be a smooth projective curve 
of genus $g \geq 2$ defined over $K$.
Let $J$ be the Jacobian of $C$ and $r$ the rank of the Mordell--Weil group $J(K)$.
Fix a rational point $P_0 \in C(K)$ and let $\jmath:C \hookrightarrow J$ be 
the Abel-Jacobi map with base point $P_0$. We use $\jmath$ to identify $C$
as a subvariety of $J$.

To explain the usual Chabauty method it is convenient to assume that
$K=\Q$. Choose a finite prime $p$.
Inside $J(\Q_p)$ it is clear that
\[
C (\Q) \subseteq C(\Q_p) \cap J(\Q) \subseteq C(\Q_p) \cap \overline{J(\Q)},
\]
where $\overline{J(\Q)}$ is the closure of $J(\Q)$ in the $p$-adic topology. 
Now $J(\Q_p)$ is a $\Q_p$-Lie group of dimension $g$, and
$\overline{J(\Q)}$ is a $\Q_p$-Lie subgroup of dimension at most $r$. 
Moreover, $C(\Q_p)$ is a $1$-dimensional submanifold of $J(\Q_p)$.
If $r+1 \leq g$ then we expect that the intersection 
$C(\Q_p) \cap \overline{J(\Q)}$ is finite. It turns out
that this intersection is indeed finite if $r \leq g-1$
and Coleman \cite{Co} gives a bound for the cardinality of this
intersection under some further (but mild)
 hypotheses. Moreover, in practice
this intersection can be computed to any desired $p$-adic accuracy. 

Now we return to the general setting by letting $K$ be a number field of degree $d$.
Define the Weil restrictions
\begin{equation}\label{eqn:wr}
V=\Res_{K/\Q} C, \qquad A=\Res_{K/\Q} J.
\end{equation}
Then $V$ is a variety
of dimension $d$ and $A$ an abelian variety of dimension $gd$, both defined over $\Q$. 
Moreover $\jmath: C \hookrightarrow J$ descends to a morphism $V \hookrightarrow A$
defined over $\Q$ which we use to identify $V$ as a subvariety of $A$. 
The Weil restriction defines a bijection between $C(K)$ and $V(\Q)$,
and 
\[
\rank A(\Q)=\rank J(K)=r.
\]
Mimicking the previous argument
\[
 V(\Q) \subseteq  V(\Q_p) \cap \overline{A(\Q)}.
\]
Now $\overline{A(\Q)}$ is at most $r$-dimensional, $V(\Q_p)$
is $d$-dimensional and the intersection is taking place
in the $\Q_p$-Lie group $A(\Q_p)$ of dimension $gd$. 
If $r+d \leq gd$ we expect that the intersection is finite.
As Wetherell points out, this is not always true. For example,
let $C$ be a curve defined over $\Q$ with Mordell--Weil rank $\geq g$.
One normally expects that $\overline{J(\Q)}$ is $g$-dimensional. Assume
that this is the case. Then the intersection
\[
C(\Q_p) \cap \overline{J(\Q)}
\]
will contain a neighbourhood in $C(\Q_p)$ of the base point $P_0$ and so will be infinite.
Now let $V$ and $A$ be obtained from $C$ and $J$ by first base extending to number
field $K$ and then taking Weil restriction back to $\Q$. One has
a natural injection
\[
C(\Q_p) \cap \overline{J(\Q)} \hookrightarrow V(\Q_p) \cap \overline{A(\Q)}
\]
proving that the latter intersection is infinite. This is true regardless
of whether the inequality $r \leq d(g-1)$ is satisfied.
However, for a random curve $C$ defined over a number field $K$,
on the basis for the above heuristic argument, 
we expect the intersection $V(\Q_p) \cap \overline{A(\Q)}$ 
to be finite when the inequality $r \leq d(g-1)$ is satisfied.

A {\bf possibly correct} statement is the following. Let $C$ be a smooth
projective curve of genus $g \geq 2$ over a number field $K$ of degree $d$.
Suppose that for every subfield $L \subseteq K$ and for every smooth projective curve $D$
defined over $L$ satisfying $D \times_L K \cong_K C$, the inequality
\[
\rank J_D(L) \leq [L:\Q] (g-1)
\]
holds, where $J_D$ denotes the Jacobian of $D$.
Let $V$ and $A$ be given by \eqref{eqn:wr}.
Then $V(\Q_p) \cap \overline{A(\Q)}$ is finite.

\section{Preliminaries} \label{sec:prelim}
In this section we summarise various results on $p$-adic
integration that we need. The definitions and proofs can
be found in \cite{CoAb} and \cite{Colmez}. 
For an introduction to the ideas involved
in Chabauty's method
we warmly recommend
Wetherell's thesis \cite{We} and the survey paper
of McCallum and Poonen \cite{MP}, as well as Coleman's
paper~\cite{Co}.

\subsection{Integration}
Let $p$ be a (finite) rational prime. Let $K_\vv$ be a finite
extension of $\Q_p$ and $\OO_\vv$ be the ring of integers
in $K_\vv$.
Let $\cW$ be a smooth, proper
connected scheme of finite type over $\OO_\vv$
and write $W$ for the generic fibre. 
In \cite[Section II]{CoAb} Coleman describes 
how to integrate \lq\lq differentials of the second
kind\rq\rq\ on $W$. We shall however only be concerned
with global $1$-forms (i.e. differentials of the first
kind) and so shall restrict our attention
to these. Among the properties  of integration 
(see \cite[Section II]{CoAb}) 
we shall need are the following:
\begin{enumerate}
\item[(i)] $\displaystyle \int_{P}^Q \omega= -\int_Q^P \omega$,
\item[(ii)] $\displaystyle \int_Q^P \omega+\int_P^R \omega=
\int_Q^R \omega$,
\item[(iii)] $ \displaystyle \int_{Q}^{P} \omega+\omega^\prime
= \int_{Q}^{P} \omega+
\int_{Q}^{P} \omega^\prime$, 
\item[(iv)] $ \displaystyle \int_Q^P \alpha \omega= \alpha \int_Q^P \omega$, 
\end{enumerate}
for $P$, $Q$, $R \in W(K_\vv)$, 
global $1$-forms $\omega$, $\omega^\prime$
on $W$, and $\alpha \in K_\vv$.
We shall also need the \lq\lq change of variables formula\rq\rq\ 
\cite[Theorem 2.7]{CoAb}:
if $\cW_1$, $\cW_2$ are smooth, proper connected schemes of 
finite type over $\OO_\vv$ 
and $\varrho:W_1 \rightarrow W_2$
is a morphism of their generic fibres then
\[
\int_{Q}^P \varrho^* \omega= \int_{\varrho(Q)}^{\varrho(P)} \omega
\]
for all global $1$-forms $\omega$ on $W_2$ and $P$, $Q \in W_1(K_\vv)$.

Now let $A$ be an abelian variety of dimension $g$ over $K_\vv$,
and write $\Omega_A$ for the $K_\vv$-space of global 
$1$-forms on $A$. Consider the pairing
\begin{equation}\label{eqn:abpairing}
\Omega_A \times A(K_\vv) \rightarrow K_\vv, \qquad
(\omega,P) \mapsto \int_0^P \omega.
\end{equation}
This pairing is bilinear. It is $K_\vv$-linear on the left
by (iii) and (iv). It is $\Z$-linear 
on the right; this is a straightforward
consequence \cite[Theorem 2.8]{CoAb} of the
\lq\lq change of variables formula\rq\rq. The kernel on the 
left is $0$ and on the right is the torsion subgroup
of $A(K_\vv)$; see \cite[III.7.6]{Bourbaki}. 

\subsection{Notation}
Henceforth we shall be concerned with curves over number fields
and their Jacobians. We fix once and for all the following notation:

\medskip
\begin{tabular}{lll}
$K$ & $\qquad$ & {a number field,}\\
$C$ &  & {a smooth projective absolutely irreducible curve} \\
& & {defined over $K$, of genus $\geq 2$,}\\
$J$ &  & {the Jacobian of $C$,} \\
$\vv$ & & {a non-archimedean prime of $K$, of good reduction for $C$,}\\
$K_\vv$ & & {the completion of~$K$ at~$\vv$,}\\
$k_\vv$ & & {the residue field of $K$ at $\vv$,}\\
$\OO_\vv$ & & {the ring of integers in $K_\vv$,}\\
$x \mapsto \tilde{x}$ && {the natural map $\OO_\vv\rightarrow k_\vv$,}\\
$\cC_\vv$ & & {a minimal regular proper model for $C$ over $\OO_\vv$,}\\
$\tilde{C}_\vv$ & & {the special fibre of $\cC_\vv$ at $\vv$,}\\
$\Omega_{C/K_\vv}$ & & {the $K_\vv$-vector space of global 
$1$-forms on $C$.}
\end{tabular}

\subsection{Integration on Curves and Jacobians}
For any field extension $M/K$ (not necessarily finite),
we shall write $\Omega_{C/M}$
and $\Omega_{J/M}$ for the $M$-vector spaces of 
global $1$-forms on $C/M$ and $J/M$ respectively.
We shall assume the existence of some $P_0 \in C(K)$.
Corresponding to $P_0$ 
is the Abel--Jacobi map, 
\[
\jmath: C \hookrightarrow J, \qquad P \mapsto [P-P_0].
\]
It is well-known that the pull-back 
$\jmath^*:\Omega_{J/K} \rightarrow \Omega_{C/K}$
is an isomorphism of $K$-vector spaces 
\cite[Proposition 2.2]{Milne}. Moreover
any two Abel--Jacobi maps differ by a translation on $J$.
As $1$-forms on $J$ are translation invariant,
the map $\jmath^*$ is independent of the choice of $P_0$
(see \cite[Section 1.4]{We}). 
Let $\vv$ be a non-archimedean place for $K$.
The isomorphism $\jmath^*$ extends to an isomorphism
$\Omega_{J/K_\vv} \rightarrow \Omega_{C/K_\vv}$,
which we shall also denote by $\jmath^*$.
For any global $1$-form $\omega \in \Omega_{J/K_\vv}$
and any two points $P$, $Q \in C(K_\vv)$ we have
\[
\int_Q^P \jmath^* \omega=\int_{\jmath{Q}}^{\jmath{P}} \omega
=\int_0^{[P-Q]} \omega,
\] 
using the properties of integration above.
We shall henceforth use $\jmath^*$ to identify $\Omega_{C/K_\vv}$ with
$\Omega_{J/K_\vv}$. With this identification, the pairing
\eqref{eqn:abpairing} with $J=A$ gives the bilinear pairing
\begin{equation}\label{eqn:pairing}
\Omega_{C/K_\vv} \times J(K_\vv) \rightarrow K_\vv, \qquad
\left(\, \omega,\; \left[\sum P_i -  Q_i\right] \right)
\longmapsto
\sum \int_{Q_i}^{P_i} \omega,
\end{equation}
whose kernel on the right is $0$ and on the left is the
torsion subgroup of $J(K_\vv)$.
We ease notation a little by defining, for divisor class 
$D=\sum P_i-Q_i$ of
degree
$0$, the integral
\[
\int_D \omega=\sum \int_{Q_i}^{P_i} \omega.
\]
Note that this integral depends on the equivalence class of $D$ 
and not on its  
decomposition as $D=\sum P_i-Q_i$.

\subsection{Uniformizers}\label{sec:uniformizers}
The usual Chabauty approach when studying rational points in
a residue class is to work with a  local coordinate (defined
shortly) and create power-series equations in terms of the local
coordinate whose solutions, 
roughly speaking, contain the rational points.
In our situation we find it more convenient to 
shift the local coordinate 
so that it becomes a uniformizer at 
a rational point in the residue class.
Fix a non-archimedean place $\vv$ of good reduction for $C$,
and a minimal regular proper model $\cC_\vv$ for $C$ over $\vv$. 
Since our objective is explicit computation,
we point out that in our case of good reduction, such a model is  
simply a system of
equations for the non-singular curve that reduces to a non-singular
system modulo $\vv$.
Let $Q \in C(K)$ and let $\tilde{Q}$ be its reduction on
the special fibre $\tilde{C}_\vv$. 
Choose a rational function $s_Q \in K(C)$ so that 
the maximal ideal in $\OO_{\cC_\vv,\tilde{Q}}$ is $(s_Q,\pi)$,
where $\pi$ is a uniformizing element for $K_\vv$. 
The function $s_Q$
is called  \cite[Section 1]{LT} a {\em local coordinate} at $Q$.
Let $t_Q=s_Q-s_Q(Q)$. We shall refer to $t_Q$,
constructed as above, as a {\em well-behaved uniformizer} 
at $Q$. The reason for the name will be clear from Lemma~\ref{lem:wb}
below. 

Before stating the lemma we define the {\em $\vv$-unit ball
around $Q$} to be
\begin{equation}\label{eqn:vub}
\cB_\vv(Q)=\{P \in C(K_\vv): \tilde{P}=\tilde{Q}\}.
\end{equation}
\begin{lem}\label{lem:wb}
\begin{enumerate}
\item[(i)] $t_Q$ is a uniformizer at $Q$,
\item[(ii)] $\tilde{t}_Q$ is a uniformizer at $\tilde{Q}$,
\item[(iii)] 
Let 
$\pi$ be a uniformizing element for $K_\vv$.  Then $t_Q$ 
is regular and injective on $\cB_\vv(Q)$.
Indeed, $t_Q$ defines a bijection  
between $\cB_\vv(Q)$
and $\pi \OO_\vv$, given by $P \mapsto t_Q(P)$. 
\end{enumerate}
\end{lem}
\begin{proof}
Parts (i) and (ii) are clear from the construction. Part (iii)
is standard; see for example \cite[Section 1]{LT}
or \cite[Sections 1.7, 1.8]{We}.
\end{proof}

\subsection{Estimating Integrals on Curves}\label{subsec:eval}
\begin{lem}\label{lem:approx} 
Let $p$ be an odd rational prime that does not ramify in $K$. 
Let $\vv$ be a place of $K$ above $p$.
Fix a minimal regular model $\cC_\vv$ for $C$ over $\OO_\vv$.
Let  $Q \in C(K_\vv)$ and let $t_Q \in K(C)$ be a 
well-behaved uniformizer at $Q$. 
Let $\omega \in \Omega_{\cC_\vv/\OO_\vv}$, and write
\begin{equation}\label{eqn:alpha}
\alpha=\frac{\omega}{dt_Q}\Bigr\rvert_{t_Q=0}.
\end{equation}
Then $\alpha \in \OO_\vv$. Moreover,
for all $P \in \cB_\vv(Q)$, 
\begin{equation}\label{eqn:approx}
\int_Q^P \omega = \alpha \cdot t_Q(P)  + \beta \cdot t_Q(P)^2
\end{equation}
for some $\beta \in \OO_\vv$ (which depends on $P$).
\end{lem}
\begin{proof}
We can expand $\omega$ (after viewing it
as an element in $\Omega_{\hat{\OO}_Q}$)
as a formal power series
\begin{equation}\label{eqn:exp}
\omega=(\alpha_0+ \alpha_1 t_Q + \alpha_2 t_Q^2+ \cdots) dt_Q,
\end{equation}
where the coefficients $\alpha_i$ are all integers in $K_\vv$
(see for example \cite[Proposition 1.6]{LT} or 
\cite[Chapters 1.7, 1.8]{We}); here we have not
used the assumption that $t_Q(Q)=0$, merely that $t_Q$
is a local coordinate at $Q$. 
We note that $\alpha=\alpha_0$ and hence integral.

Let $P \in \cB_\vv(Q)$.
We can now evaluate the integral
(see for example \cite[Proposition 1.3]{LT})
\[
\int_Q^P \omega=\sum_{j=0}^\infty \frac{\alpha_j}{j+1} t_Q(P)^{j+1};
\]
the infinite series converges since $\ord_\vv(t_Q(P)) \geq 1$
by part (iii) of Lemma~\ref{lem:wb}. 
Thus \eqref{eqn:approx} holds with
\[
\beta=\sum_{j=0}^\infty \frac{\alpha_{j+1}}{j+2} t_Q(P)^{j}.
\]
To complete the proof we must show that $\beta$ is integral.
Thus, it is sufficient to show that
\[
\ord_\vv (j+2) \leq  j
\] 
for all $j \geq 0$. But $K_\vv/\Q_p$ is unramified, and so
$\ord_\vv (j+2)=\ord_p(j+2)$. Hence we need to
show that $\ord_p(j+2) \leq j$ for all $j \geq 0$ and
all odd primes $p$. This is now an easy exercise.
\end{proof}

\section{Chabauty in a Single Unit Ball}\label{sec:chab}

Let $C$ be a smooth projective curve over a number field $K$. 
Let $J$ be the Jacobian of $C$ and write $r$ for the
rank of the Mordell--Weil group $J(K)$. Let
$D_1,\cdots,D_r$ be a basis for a free subgroup of finite
index in $J(K)$. 
Let $p$ a rational prime satisfying the following: \label{assumptions}
\begin{enumerate}
\item[(p1)] $p$ is odd
\item[(p2)] $p$ is unramified in $K$,
\item[(p3)] every prime $\vv$ of $K$ above $p$ is
a prime of good reduction for the curve $C$.
\end{enumerate}

For each $\vv \mid p$ we fix once and for all a minimal regular proper model
$\cC_\vv$ for $C$ over $\OO_\vv$. 
Let $Q \in C(K)$. 
For $\vv \mid p$, let $\cB_\vv(Q)$ be as in \eqref{eqn:vub}, and define the {\em $p$-unit
ball around $Q$} to be
\begin{equation}\label{eqn:pub}
\cB_p(Q)=\prod_{\vv \mid p} \cB_\vv(Q).
\end{equation} 
We will shortly give
a criterion for a point $Q \in C(K)$ to be the unique $K$-rational point
in its own $p$-unit ball.

To state our criterion---Theorem~\ref{thm:1} below---we 
need to define a pair of matrices $T$ and $A$. 
The matrix $T$ depends on the basis $D_1,\dotsc,D_r$. The matrix $A$
depends on the point $Q \in C(K)$.
Let $\vv_1,\dotsc,\vv_n$ be the places of $K$ above $p$.
For each place $\vv$ above $p$ we fix once and for all a $\Z_p$-basis
$\theta_{\vv,1},\dots,\theta_{\vv,d_\vv}$ for $\OO_\vv$, where 
$d_\vv=[K_\vv:\Q_p]$. Of course $d_\vv=[\OO_\vv:\Z_p]=[k_\vv:\F_p]$ as $p$ is unramified
in $K$. We also  choose 
an $\OO_\vv$-basis
$\omega_{\vv,1},\dotsc,\omega_{\vv,g}$ for $\Omega_{\cC_\vv/\OO_\vv}$.

Now fix $\vv$ above $p$, and let $\omega \in \Omega_{\cC_\vv/\OO_\vv}$.
Let
\begin{equation}\label{eqn:tau}
\tau_j=\int_{D_j} \omega, \qquad j=1,\dotsc,r.
\end{equation}
Write
\begin{equation}\label{eqn:tij}
\tau_j=\sum_{i=1}^{d_\vv} t_{ij} \theta_{\vv,i}, \qquad t_{ij} \in \Q_p.
\end{equation}
Let 
\begin{equation}\label{eqn:Tvomega}
T_{\vv,\omega}=(t_{ij})_{i=1,\dotsc,d_\vv,\; j=1,\dotsc,r};
\end{equation}
that is $T_{\vv,\omega}$ is the $d_\vv \times r$ matrix with entries $t_{ij}$.
Recall that $\omega_{\vv,1},\dotsc,\omega_{\vv,g}$ is a basis for $\Omega_{\cC_\vv/\OO_\vv}$.
Let
\begin{equation}\label{eqn:Tv}
T_{\vv}=
\begin{pmatrix}
T_{\vv,\omega_{\vv,1}} \\ 
T_{\vv,\omega_{\vv,2}} \\
\vdots \\ 
T_{\vv,\omega_{\vv,g}}
\end{pmatrix}; 
\end{equation}
this is a $g d_\vv \times r$ matrix with entries in $\Q_p$. We now define
the matrix $T$ needed for our criterion below:
\begin{equation}\label{eqn:T}
T=
\begin{pmatrix}
T_{\vv_1} \\
T_{\vv_2} \\
\vdots \\
T_{\vv_n}
\end{pmatrix}.
\end{equation}
Note that $T$ is $gd \times r$ matrix with entries in $\Q_p$ where
$d=[K:\Q]=d_{\vv_1}+\cdots+d_{\vv_n}$.

Let $Q \in C(K)$. We now define the second matrix $A$ needed
to state our criterion for $C(K) \cap \cB_p(Q)=\{ Q\}$. 
For each place $\vv$ of $K$ above $p$, we have chosen a minimal proper
regular model $\cC_\vv$. Let $t_Q$ be a well-behaved uniformizer at $Q$
as defined in Subsection~\ref{sec:uniformizers}.
Let $\omega \in \Omega_{\cC_\vv/\OO_\vv}$
and let 
$\alpha$ be given by~\eqref{eqn:alpha}.
By Lemma~\ref{lem:approx}, $\alpha \in \OO_\vv$. Recall that we have
fixed a basis $\theta_{\vv,1},\dotsc,\theta_{\vv,d_\vv}$ for $\OO_\vv/\Z_p$.
Write 
\begin{equation}\label{eqn:aij}
\alpha \cdot \theta_{\vv,j}=\sum_{i=1}^{d_\vv} a_{ij} \theta_{\vv,i},
\qquad j=1,\dots,d_\vv,
\end{equation}
 with $a_{ij} \in \Z_p$.
Let 
\begin{equation}\label{eqn:Avomega}
A_{\vv,\omega}=(a_{ij})_{i,j=1,\dotsc,d_\vv}.
\end{equation}
Let
\begin{equation}\label{eqn:Av}
A_\vv=\begin{pmatrix}
A_{\vv,\omega_1} \\
A_{\vv,\omega_2} \\
\vdots \\
A_{\vv,\omega_g}
\end{pmatrix};
\end{equation}
this is a $d_\vv g \times d_\vv$ matrix with entries in $\Z_p$. Let
\begin{equation}\label{eqn:AQ}
A=\left(
\begin{array}{cccc}
A_{\vv_1} & \mathbf{0} &  \hdots & \mathbf{0} \\
\mathbf{0} & A_{\vv_2} &  \hdots & \mathbf{0} \\
\vdots & \vdots & \ddots  & \vdots \\
\mathbf{0} & \mathbf{0} &  \hdots & A_{\vv_n}
\end{array}
\right).
\end{equation}
Then $A$ is a $dg \times d$ matrix with entries in $\Z_p$.

We use $A$ and $T$ to define a matrix $M_p(Q)$ in terms
of which we will express our criterion for $Q$ to be the unique
rational point in its $p$-unit ball. Choose a non-negative integer $a$
such that $p^a T$ has entries in $\Z_p$. Let $U$ be a unimodular matrix
with entries in $\Z_p$ such that $U (p^a T) $ is in Hermite Normal Form. 
\label{HNF}
Let $h$ be the number of zero rows of $U (p^a T) $. Let
$M_p(Q)$ be the $h \times d$ matrix (with entries in $\Z_p$)
formed by the last $h$ rows of $UA$.
\begin{thm}\label{thm:1}
With the above assumptions and notation, 
denote
by $\tilde{M}_p(Q)$ the matrix with entries in $\F_p$ obtained
by reducing $M_p(Q)$ modulo $p$. If $\tilde{M}_p(Q)$
has rank $d$ then $C(K) \cap \cB_p(Q) =\{Q\}$.
\end{thm}

\noindent \textbf{Remarks.} 

\begin{enumerate}
\item[(i)] Let $\mathbf{u}_1,\dots,\mathbf{u}_h$ be a $\Z_p$-basis
for the kernel of the homomorphism of $\Z_p$-modules 
$\Z_p^{gd} \rightarrow \Z_p^r$ given by $p^a T$. Then
$\mathbf{u}_1 A,\dots,\mathbf{u}_h A$ span the same
$\Z_p$-module the rows of $M_p(Q)$, showing that the rank
of $\tilde{M}_p(Q)$ is independent of the choice of $U$.

\item[(ii)] Since the matrix $T$ is $gd \times r$, it is evident
that $h \geq \max(gd-r,0)$ and, very likely, $h=\max(gd-r,0)$.
Now the matrix $\tilde{M}_p(Q)$ is
$h \times d$ and so a necessary condition for the criterion to hold is that $h \geq d$. 
Thus it is sensible to apply the theorem when $gd-r \geq d$, or
equivalently when $r \leq d(g-1)$. 
\item[(iii)] In practice, we do not 
compute the matrix $T$ exactly, merely an approximation to it. Thus
we will not be able to provably determine $h$ unless $h=\max(gd-r,0)$.
\end{enumerate}

\begin{proof}[Proof of Theorem~\ref{thm:1}]
Suppose
that $P \in C(K) \cap \cB_p(Q)$.
We need to show that $P=Q$.

Let $m$ be the index
\begin{equation}\label{eqn:index}
m:= [J(K): \langle D_1,\cdots,D_r \rangle].
\end{equation}
There are integers $n_1^\prime,\dotsc,n_r^\prime$ such that
\begin{equation}\label{eqn:mw}
m(P-Q)=n_1^\prime D_1+\cdots+n_r^\prime D_r,
\end{equation}
where the equality takes place in $\Pic^0(C)$.

Let $\vv$ be one of the places $\vv_1,\dotsc,\vv_n$ above $p$. 
Recall that we have chosen a well-behaved uniformizer $t_Q$ at $Q$.
Write $z=t_Q(P)$.   
By part (iii) of Lemma~\ref{lem:wb},
 $\ord_\vv(z) \geq 1$. We will show that
$z=0$, and so again by part (iii) of Lemma~\ref{lem:wb}, $P=Q$ which 
is what we want
to prove.

We write 
\[
z=z_{\vv,1}\theta_{\vv,1}+\cdots+z_{\vv,d_\vv} \theta_{\vv,d_\vv},
\]
where $z_{\vv,i} \in \Z_p$. 
As $\tilde{\theta}_{\vv,1},\dotsc,\tilde{\theta}_{\vv,d_\vv}$
is a basis for $k_\vv/\F_p$ and $\ord_\vv(z) \geq 1$, 
we see that $\ord_\vv(z_{\vv,i}) \geq 1$ for $i=1,\dotsc,d_\vv$.
Let 
\begin{equation}\label{eqn:sv}
s_\vv=\min_{1 \leq i \leq d_\vv}{\ord_p(z_{\vv,i})}.
\end{equation}
We will show that $s_\vv=\infty$, which implies that $z_{\vv,i}=0$ for $i=1,\dotsc,d_\vv$
and so $z=0$ as required.
For now we note that $s_\vv \geq 1$.

Now fix an $\omega \in \Omega_{\cC_\vv/\OO_\vv}$ and let 
$\alpha \in \OO_\vv$ be as in 
Lemma~\ref{lem:approx}; by that lemma
\[
\int_Q^P \omega = \alpha  z +\beta z^2,
\]
for some $\beta \in \OO_\vv$.
However, by equation~\eqref{eqn:mw} and the properties of integration
explained in Section~\ref{sec:prelim},
\[
m\int_Q^P \omega=n_1^\prime \tau_1+\cdots+ n_r^\prime \tau_r,
\]
where the $\tau_j$ are given in \eqref{eqn:tau}.
Let $n_i=n_i^\prime/m \in \Q$. 
Thus
\[
\int_Q^P \omega=n_1 \tau_1+\cdots + n_r \tau_r. 
\]
Hence
\begin{equation}\label{eqn:linquad}
n_1 \tau_1+\cdots+ n_r \tau_r=
\alpha (z_{\vv,1}\theta_{\vv,1}+\cdots+z_{\vv,d_\vv} \theta_{\vv,d_\vv})
+\beta (z_{\vv,1}\theta_{\vv,1}+\cdots+z_{\vv,d_\vv} \theta_{\vv,d_\vv})^2.
\end{equation}
From \eqref{eqn:linquad} and \eqref{eqn:sv}
we obtain
\begin{equation}\label{eqn:lin}
n_1 \tau_1+\cdots+n_r \tau_r \equiv 
z_{\vv,1}(\alpha \theta_{\vv,1})+\cdots+z_{\vv,d_\vv} (\alpha \theta_{\vv,d_\vv}) \pmod{p^{2 s_\vv}\OO_\vv}.
\end{equation}
Write
\[
\mathbf{n}=\begin{pmatrix}
n_1 \\ n_2 \\ \vdots \\ n_r
\end{pmatrix}, \qquad
\mathbf{z}_\vv=
\begin{pmatrix}
z_{\vv,1} \\ z_{\vv,2} \\ \vdots \\ z_{\vv,d_\vv}
\end{pmatrix},
\]
and note that the entries of $\mathbf{n}$ are in $\Q$,
and the entries of $\mathbf{z}_\vv$ are in $p^{s_\vv}\Z_p$.
Recall that we have expressed
 $\tau_j=\sum t_{ij} \theta_{\vv,i}$ in \eqref{eqn:tij}
and $\alpha \cdot \theta_{\vv,j}=\sum a_{ij} \theta_{\vv,i}$
in \eqref{eqn:aij}
where $t_{ij}$ are in $\Q_p$ and the $a_{ij}$ are in $\Z_p$. 
Substituting in \eqref{eqn:lin} and comparing the coefficients for $\theta_{\vv,i}$ we obtain
\[
T_{\vv,\omega} \mathbf{n}
\equiv
A_{\vv,\omega} 
\mathbf{z}_\vv
\pmod{p^{2 s_\vv}},
\]
where $T_{\vv,\omega}$ 
and $A_{\vv,\omega}$ are respectively given in \eqref{eqn:Tvomega}
and \eqref{eqn:Avomega}.

Let $T_\vv$ and $A_\vv$ be as given
in \eqref{eqn:Tv} and \eqref{eqn:Av} respectively. Then
\[
T_{\vv} \mathbf{n}
\equiv
A_{\vv} \mathbf{z}_\vv
\pmod{p^{2 s_\vv}}.
\]

Now let
\[
\mathbf{z}=\begin{pmatrix} \mathbf{z}_{\vv_1} \\ 
\mathbf{z}_{\vv_2}\\
\vdots \\ \mathbf{z}_{\vv_n}
\end{pmatrix}.
\]
Then $\mathbf{z}$ is of length $d=[K:\Q]$ with entries in $p\Z_p$.
Write
\begin{equation}\label{eqn:s}
s=\min_{\vv=\vv_1,\dots,\vv_n}{s_{\vv}} = \min_{i,j}\; \ord_{\vv_j}(z_{i,\vv_j}),
\end{equation}
where the $s_\vv$ are defined in \eqref{eqn:sv}.
Clearly $s \geq 1$. It is sufficient to show that $s=\infty$ since
then all of the $z_{i,\vv_j}=0$ implying that $P=Q$.

Let $T$ and $A$ be as given in \eqref{eqn:T} and \eqref{eqn:AQ}.
Then
\begin{equation}\label{eqn:Tcong}
T \mathbf{n}  \equiv A \mathbf{z} \pmod{p^{2s}},
\end{equation}
where we note once again that $T$ is $dg\times r$ with entries in $\Q_p$
and $A$ is $dg \times d$ with entries in $\Z_p$.

Let $U$, $M_p(Q)$, $h$ be as in the paragraph preceding the statement of the
theorem. 
Suppose that $\tilde{M}_p(Q)$ has rank $d$. 
Suppose $s< \infty$ and we will derive a contradiction. Recall that the
last $h$ rows of $UT$ are zero. From 
 \eqref{eqn:Tcong} we have that
$M_p(Q) \mathbf{z} \equiv 0 \pmod{p^{2s}}$
since, by definition, $M_p(Q)$ is the matrix formed by the
last $h$ rows of $UA$. In particular $M_p$ has coefficients in $\Z_p$ since
both $U$ and $A$ have coefficients in $\Z_p$. From the definition of $s$ in 
\eqref{eqn:s} we can write $\mathbf{z}=p^s \mathbf{w}$ where the entries
of $\mathbf{w}$ are in $\Z_p$ and $\mathbf{w} \not \equiv \mathbf{0} \pmod{p}$.
However,
$M_p(Q) \mathbf{w} \equiv 0 \pmod{p^s}$,
and as $s \geq 1$, we have that $M_p(Q) \mathbf{w} \equiv 0 \pmod{p}$. Since 
$\mathbf{w} \in \Z_p^d$, if $\tilde{M}_p(Q)$ has rank $d$ then
$\mathbf{w} \equiv \mathbf{0} \pmod{p}$, giving the desired contradiction.
\end{proof}

\section{Chabauty and the Mordell--Weil Sieve}\label{sec:MW}
Theorem~\ref{thm:1} gives a criterion for showing that 
a given $K$-rational point $Q$ on $C$ is the unique $K$-rational
point in its $p$-unit ball, for a rational prime $p$
satisfying certain conditions.

Let $L_0$ be a subgroup of $J(K)$ of finite index containing
the free subgroup $L$ generated by $D_1,\dotsc,D_r$ of the previous
section. We can take $L_0=L$ but for our purpose it is preferable to 
include the torsion subgroup of $J(K)$ in $L_0$. 
The usual $p$-saturation method \cite{Si},
\cite{Sikphd}, 
\cite{FS} shows how to enlarge $L_0$ so that its
index in $J(K)$ is not divisible by any given small prime $p$.
One expects, after checking $p$-saturation for all small primes
$p$ up to some bound that $L_0$ is in fact equal to $J(K)$.
However, proving that $J(K)=L_0$ requires an explicit theory
of heights on the Jacobian $J$. This is not yet available
for Jacobians of curves of genus $\geq 3$. For Jacobians of 
curves of genus $2$ there is an explicit theory of heights
\cite{Flynnh}, \cite{FS}, \cite{Stollh1}, \cite{Stollh2},
though the bounds over any  number fields other than the rationals
are likely to be impractically large.

As usual we assume the existence of some $P_0 \in C(K)$
and denote the associated Abel--Jacobi map by $\jmath$.
It is easiest to search for $K$-rational points on $C$ 
by taking small linear combinations of the generators of $L_0$
and checking if they are in the image of $\jmath$.
In this way we determine a set $\cK \subseteq C(K)$
of known $K$-rational points, and the challenge is to
show that $\cK$ is in fact equal to $C(K)$. In this section
we show how to combine our Theorem~\ref{thm:1} with
an adaptation of the Mordell--Weil sieve to give
a practical criterion for $\cK$ to be equal to $C(K)$.
The usual Mordell--Weil sieve \cite{BrE}, \cite{BS}, \cite{BS3}, 
\cite{integral}
assumes knowledge of the
full Mordell--Weil group. Our adaptation takes careful account
of the fact that we are working with a subgroup of finite (though unknown)
index.

Before we give the details we point 
that substantial improvements can be made to
 the version of the Mordell--Weil sieve 
outlined below. It has 
certainly been sufficient for the examples we have 
computed so far (including the ones detailed in the next
section). But we expect that for some other examples it will be necessary
(though not difficult)
to incorporate the improvements to the Mordell--Weil sieve
found in the papers
of Bruin and Stoll \cite{BS}, \cite{BS3}. 

Before stating our criterion we need to set up some notation.
Here it is convenient to use the same symbol to denote several
different Abel-Jacobi maps associated to the fixed $K$-rational
point $P_0$ and its images on special fibres.
Let  $\vv$ be of good reduction for $C$. 
Denote by $\red$ the natural maps 
\[
\red: C(K) \rightarrow C(k_\vv), \qquad
\red: J(K) \rightarrow J(k_\vv).
\]
 Denote by
$\jmath$ the Abel--Jacobi map $C(k_\vv) \rightarrow J(k_\vv)$
associated to $\tilde{P}_0$. 

\begin{lem}\label{lem:mw}
Let $L_0$ be a subgroup of $J(K)$ of finite index $n=[J(K):L_0]$. 
Let $P_0 \in C(K)$ and let $\jmath$ denote 
the Abel-Jacobi maps associated to $P_0$ as above. 
Let 
$\vv_1,\dotsc,\vv_s$ be places of $K$,
 such that each $\vv=\vv_i$ satisfies the  
two conditions:
\begin{enumerate}
\item[($\upsilon 1$)] $\vv$ is a place of good reduction for $C$, 
\item[($\upsilon 2$)] the index $n$ is coprime to $\# J(k_\vv)$.
\end{enumerate}
To ease notation, write $k_i$ for the residue field $k_{\vv_i}$.
Define inductively a sequence of subgroups
\[
L_0 \supseteq L_1 \supseteq L_2 \supseteq L_3 \supseteq \dots \supseteq L_s
\] 
and finite subsets $W_0,W_1,\dotsc,W_s \subset L_0$ as follows. Let
$W_0=\{\mathbf{0}\}$. Suppose we have defined
$L_i$ and $W_i$ where $i \leq s-1$.
Let $L_{i+1}$ be the kernel
of the composition
\[
L_i \hookrightarrow J(K) \rightarrow  J(k_{i+1}).
\]
To define $W_{i+1}$
choose a complete set $\cQ$ of coset representatives for $L_{i}/L_{i+1}$,
and let
\[
W_{i+1}^\prime=\{\mathbf{w}+\mathbf{q}: \mathbf{w} \in W_i \quad \text{and} \quad 
\mathbf{q} \in \cQ\}. 
\]
Let
\[
W_{i+1}=\{ \mathbf{w} \in W_{i+1}^\prime : \red(\mathbf{w}) \in \jmath(C(k_{i+1}))\}.
\]
Then for every $i=0,\dotsc,s$, and every $Q\in C(K)$, there is some $\mathbf{w} \in W_i$
such that 
\begin{equation}\label{eqn:ind}
n(\jmath(Q)-\mathbf{w}) \in L_i.
\end{equation}
\end{lem}
\begin{proof}
The proof is by induction on $i$. Since 
$L_0$ has index $n$ in $J(K)$, \eqref{eqn:ind} is true with $\mathbf{w}=0$.
Let $i \leq s-1$. Suppose $Q \in C(K)$, $\mathbf{w}^\prime \in W_i$
and $\mathbf{l}^\prime \in L_i$ satisfy
\begin{equation}\label{eqn:inter}
n(\jmath(Q)-\mathbf{w}^\prime) = \mathbf{l}^\prime.
\end{equation}
By definition of $L_{i+1}$, the quotient group $L_i/L_{i+1}$
is isomorphic to a subgroup of $J(k_{i+1})$. It follows from
 assumption ($\vv 2$) that $n$ is coprime to the order of $L_i/L_{i+1}$.
Recall that $\cQ$ was defined as a complete set of coset
representatives for $L_i/L_{i+1}$.
Thus $n \cQ$ is also a set of coset representatives. Hence
we may express $\mathbf{l}^\prime \in L_i$ as
\[
\mathbf{l}^\prime =n \mathbf{q}+ \mathbf{l}
\]
where $\mathbf{q} \in \cQ$ and $\mathbf{l} \in L_{i+1}$. Let 
$\mathbf{w}=\mathbf{w}^\prime+\mathbf{q}$. Then $\mathbf{w} \in W_{i+1}^\prime$.
By \eqref{eqn:inter},  we see that
\[
n(\jmath(Q)-\mathbf{w})=\mathbf{l}^\prime -n \mathbf{q} = \mathbf{l} \in L_{i+1}.
\]
To complete the inductive argument all we need to show is that $\mathbf{w} \in W_{i+1}$,
or equivalently that $\red(\mathbf{w}) \in \jmath (C(k_{i+1}))$. However,
since $L_{i+1}$ is contained in the kernel of $\red : J(K) \rightarrow J(k_{i+1})$,
we see that
\[
n (\jmath(\tilde{Q})-\red(\mathbf{w}))=0 \qquad \text{in $J(k_{i+1})$}.
\]
Using the fact that $n$ is coprime to $\# J(k_{i+1})$ once again gives
$\red(\mathbf{w}) =\jmath(\tilde{Q})$ as required.
\end{proof}

\begin{thm}\label{thm:mwc}
We continue with the above notation and assumptions.
Let $L_0\supseteq L_1 \supseteq \cdots \supseteq L_s$
and $W_0,\dotsc,W_s$ be the sequences constructed in Lemma~\ref{lem:mw}. 
Let $\cK$ be a subset of $C(K)$.
Let $P_0 \in \cK$ and let $\jmath$ denote 
the maps associated to $P_0$ as above. 
Suppose that for every $\mathbf{w} \in W_s$ there is a point 
$Q \in \cK$ and a prime $p$ such that  the following conditions hold:
\begin{enumerate}
\item[(a)] $p$ satisfies conditions (p1)--(p3) on page~\pageref{assumptions},
\item[(b)] (in the notation of the previous section) 
the matrix $\tilde{M}_p(Q)$ has rank $d$,
\item[(c)] the kernel of the homomorphism
\begin{equation}\label{eqn:proker}
J(K) \longrightarrow \prod_{\vv \mid p} J(k_\vv)
\end{equation}
contains both the group $L_s$ and the difference $\jmath(Q)-\mathbf{w}$,
\item[(d)] the index $n=[J(K):L_0]$ is coprime to the orders of the groups $J(k_\vv)$
for $\vv \mid p$.
\end{enumerate}
Then $C(K)=\cK$.
\end{thm}
\begin{proof}
Suppose that $P \in C(K)$. We would like to show that $P \in \cK$.
By Lemma~\ref{lem:mw}, there is some $\mathbf{w} \in W_s$ such that 
$n(\jmath(P)-\mathbf{w}) \in L_s$.
Let $Q \in \cK$ and prime $p$ satisfy conditions (a)--(d) of the theorem.
By (c), $L_s$ is contained in the kernel of \eqref{eqn:proker} and hence 
\[
n (\jmath(\tilde{P}) - \red(\mathbf{w}))=0
\] 
in $J(k_\vv)$ for all $\vv \mid p$. Since $p$ satisfies assumption (d),
it follows that
\[
\jmath(\tilde{P})-\red(\mathbf{w})=0
\]
in $J(k_\vv)$ for all $\vv \mid p$. But by assumption (c) again,
\[
\jmath(\tilde{Q})-\red(\mathbf{w})=0
\]
in $J(k_\vv)$ for all $\vv \mid p$. It follows that $\tilde{P}=\tilde{Q}$
in $C(k_\vv)$ for all $\vv \mid p$.
Hence
$P \in \cB_p(Q)$ where $\cB_p(Q)$ is the $p$-unit ball around $Q$
defined in \eqref{eqn:pub}. By assumption (b) and Theorem~\ref{thm:1}
we see that $P=Q \in \cK$ completing the proof.
\end{proof}

\section{The Generalized Fermat Equation with Signature $(2,3,10)$}\label{sec:Fermat}

Let $p$, $q$, $r \in \Z_{\geq 2}$. The equation
\begin{equation}\label{eqn:FCgen}
x^p+y^q=z^r
\end{equation}
is known as the Generalized Fermat equation (or the Fermat--Catalan equation)
 with signature $(p,q,r)$.
As in Fermat's Last Theorem, one is interested in integer solutions
$x$, $y$, $z$. Such a solution is called {\em non-trivial} if
$xyz \neq 0$, and {\em primitive} if $x$, $y$, $z$ are coprime.
Let $\chi=p^{-1}+q^{-1}+r^{-1}$. The parametrization
of non-trivial primitive solutions for $(p,q,r)$ with
$\chi \geq 1$ has now been completed \cite{Ed}.
The Generalized Fermat Conjecture \cite{Da97}, \cite{DG}
is concerned with the case $\chi<1$.
It states that the only non-trivial primitive solutions to
\eqref{eqn:FCgen} with $\chi<1$ are
\begin{gather*}
1+2^3 = 3^2, \quad 2^5+7^2 = 3^4, \quad 7^3+13^2 = 2^9, \quad
2^7+17^3 = 71^2, \\
3^5+11^4 = 122^2, \quad 17^7+76271^3 = 21063928^2, \quad
1414^3+2213459^2 = 65^7, \\
9262^3+15312283^2 = 113^7, \quad
43^8+96222^3 = 30042907^2, \quad 33^8+1549034^2 = 15613^3.
\end{gather*}
The Generalized Fermat Conjecture
has been established for many signatures $(p,q,r)$,
including for several infinite families of signatures:
Fermat's Last Theorem $(p,p,p)$ by
Wiles and Taylor \cite{W}, \cite{TW};
$(p,p,2)$ and $(p,p,3)$ by Darmon and Merel \cite{DM};
$(2,4,p)$ by Ellenberg \cite{El} and Bennett, Ellenberg and Ng
\cite{BEN};
$(2p,2p,5)$ by Bennett \cite{Ben}.
Recently, Chen and Siksek \cite{ChenS} 
have solved the Generalized Fermat equation
with signatures $(3,3,p)$ for a set of prime exponents $p$ having 
Dirichlet density $28219/44928$.
For an exhaustive survey see \cite[Chapter 14]{Cohen}. An older but
still very useful survey is \cite{Kr99}.

There is an abundance of solutions for Generalized Fermat equations
with signatures $(2,3,n)$, and so this subfamily is
particularly interesting. The condition $\chi>1$ within this subfamily coincides
with the condition $n \geq 7$. 
The cases $n=7$, $8$, $9$ are solved respectively in \cite{PSS}, \cite{Br2} 
and \cite{Br3}. The case $n=10$ appears to be the first hitherto unresolved
case within this subfamily and this of course corresponds to equation~\eqref{eqn:2310}.

In this section we solve equation~\eqref{eqn:2310}
in coprime integers $x$, $y$, $z$, thereby proving Theorem~\ref{thm:Ferm}.
We shall use the computer package {\tt MAGMA} \cite{MAGMA} for all our calculations.
In particular, {\tt MAGMA} includes implementations by 
Nils Bruin and Michael Stoll 
of $2$-descent on Jacobians of hyperelliptic curves over number fields;
the algorithm is detailed
in Stoll's paper \cite{Stoll}. {\tt MAGMA} also includes an implementation
of Chabauty for genus $2$ curves over $\Q$.

\subsection{\bf Case I: $y$ is odd}

From \eqref{eqn:2310} we immediately see that
\[
x+z^5=u^3,
\quad x-z^5=v^3, \]
where $u$, $v$ are coprime and odd.
Hence $2z^5=u^3-v^3$.

\bigskip

\noindent{\bf Case I.1: $3\nmid  z$}

Then
\[
u-v=2a^5, \qquad
u^2+uv+v^2=b^5,
\]
where $a$, $b$ are coprime integers with $z=ab$.
We now use the identity 
\begin{equation}\label{eqn:id}
(u-v)^2+3(u+v)^2=4(u^2+uv+v^2)
\end{equation}
to obtain $4a^{10}+3c^2=4b^5$, where $c=u+v$.
Dividing by $4a^{10}$, we obtain a rational point
$(X,Y)=(b/a^2,3c/2a^5)$ on the genus $2$ curve 
\[
C:\quad Y^2=3(X^5-1).
\] 
Using {\tt MAGMA} we are able to show that the Jacobian of this genus $2$
curve $C$ has Mordell-Weil rank $0$ and  torsion subgroup isomorphic to 
$\Z/2\Z$. It is immediate that $C(\Q)=\{\infty,(1,0)\}$.

Working backwards we obtain the solutions $(x,y,z)=(0,1, \pm 1)$  
to \eqref{eqn:2310}.

\bigskip

\noindent {\bf Case I.2: $3\mid z$}

Recall $2z^5=u^3-v^3$ and $u$, $v$ are odd and coprime.
Thus 
\[
u-v=2 \cdot 3^4 a^5,
\qquad      u^2+uv+v^2=3 b^5,
\] 
where $z=3ab$.
Now we use identity \eqref{eqn:id} to obtain
$4\cdot 3^8 a^{10}+3c^2=12b^5$, where $c=u+v$.
Hence we obtain a rational point $(X,Y)=(b/a^2, c/2 a^5)$ on
the genus $2$ curve
\[
C: Y^2=X^5-3^7.
\]
Let $J$ be the Jacobian of $C$.
Using {\tt MAGMA} we can show that $J(\Q)$
is free of rank $1$, with generator
\[
\left(\frac{-9+3 \sqrt{-3}}{2},\frac{81+27\sqrt{-3}}{2} \right)+
\left(\frac{-9-3 \sqrt{-3}}{2},\frac{81-27\sqrt{-3}}{2} \right)-2 \infty.
\]


Using {\tt MAGMA}'s built-in Chabauty command we find that $C(\Q)=\{\infty\}$.
Working backwards we obtain $(x,y,z)=(\pm 1,-1,0)$.

\subsection{Case II: $y$ is even}

We would now like to solve \eqref{eqn:2310} with
$y$ even, and $x$, $y$ coprime.
Replacing $x$ by $-x$ if necessary
 we obtain $x \equiv z^5 \pmod{4}$.
\[
x+z^5=2u^3, \qquad
x-z^5=4v^3,
\]
where $y=-2uv$.
Hence
\begin{equation}\label{eqn:335}
u^3-2v^3=z^5
\qquad \text{$u$, $v$ are coprime and $u$, $z$ are odd}.
\end{equation}
If $3 \mid z$ then this
equation is impossible modulo $9$. 
Hence $3 \nmid z$.

Let $\theta=\crt$.
We shall work in the number field $K=\Q(\theta)$.
This has ring of integers $\OO_K=\Z[\theta]$ with 
class number $1$.
The unit group is isomorphic to $\Z \times \Z/2\Z$
with $\epsilon=1-\theta$ a fundamental unit.

Observe that
\[
(u-v\theta)(u^2+uv \theta +v^2 \theta^2)=z^5,
\]
where the two factors on the left-hand side are coprime as 
$u$, $v$ are coprime and $z$ is neither divisible by $2$ nor $3$.
Hence 
\begin{equation}\label{eqn:uvalpha}
u-v\theta =\epsilon^s \alpha^5, \qquad
u^2+uv \theta +v^2 \theta^2=\epsilon^{-s} \beta^5,
\end{equation}
where $-2 \leq s \leq 2$ and $\alpha$, $\beta \in \Z[\theta]$ satisfy 
$z=\alpha \beta$.
We now use the identity
\[
(u-v \theta)^2+3(u+v \theta)^2=4(u^2+uv\theta+v^2 \theta^2),
\]
to obtain
\[
\epsilon^{2s} \alpha^{10}+3(u+v\theta)^2=4 \epsilon^{-s} \beta^5.
\]
Let $C_s$ be the genus $2$ curve defined over $K$ given by
\[
C_s: \quad Y^2=3(4 \epsilon^{-s} X^5- \epsilon^{2s}).
\]
We see that 
\begin{equation}\label{eqn:XY}
(X,Y)=\left(\frac{\beta}{\alpha^2},\frac{3(u+v\theta)}{\alpha^5}\right),
\end{equation}
is a $K$-rational point on $C_s$. To complete our proof of Theorem~\ref{thm:Ferm}
we need to determine $C_s(K)$ for $-2 \leq s \leq 2$. 
Let $J_s$ be the Jacobian of $C_s$. Using reduction at various places of $K$
we easily showed that the torsion subgroup of $J_s(K)$ is trivial in all cases.
The $2$-Selmer ranks of $J_s(K)$
are respectively $1$, $3$, $2$, $3$, $0$ for $s=-2,-1,0,2,1$.
We searched for $K$-rational points on each $J_s$ by
first searching for points on the associated Kummer surface.
We are fortunate to have found enough independent points in $J_s(K)$
in each case to show that the Mordell--Weil rank is equal to the 
$2$-Selmer rank. In other words we have determined a basis for a subgroup
of $J_s(K)$ of finite index, and this is given in Table~\ref{table:2310}.

\begin{table}
\begin{minipage}{\linewidth}
\caption{}
\begin{tabular}{|c|c|c|}
\hline
$s$ & basis for subgroup of $J_s(K)$ of finite index & $C_s(K)$ \\
\hline
\hline
$-2$
&
$ (\theta^2 + \theta + 1 , \theta^2 + 2\theta + 1 )-\infty $
&
$ \infty$,\\
& & 
$(\theta^2 + \theta + 1 , \pm (\theta^2 + 2\theta + 1) )$ 
\\
\hline
$-1$
&
$ (-\theta^2 - \theta - 1, 11\theta^2 + 13\theta + 17) -\infty,$
&
$\infty$,\\ 
&
$
\displaystyle 
\sum_{i=1,2}(\Phi_i, (2\theta^2 + 2\theta + 3)\Phi_i + 2\theta^2 + 3\theta + 4)
-2\infty$~\footnote{$\Phi_1$, $\Phi_2$
are the roots of $
2\Phi^2 + (\theta^2 + \theta + 2)\Phi + (\theta^2 + \theta + 2)=0$.
}, 
& 
$\displaystyle
\left( \frac{-\theta^2 - 2\theta - 1}{3} , \frac{\pm (\theta^2 - \theta + 1)}{3} \right)$, \\
&
$	
\displaystyle
\sum_{i=3,4}(\Phi_i,(4\theta^2 + 6\theta + 10)\Phi_i + 9\theta^2 + 11\theta + 13)
- 2\infty 
$~\footnote{$\Phi_3$, $\Phi_4$ are the roots of
$3\Phi^2 + (4\theta^2 + 5\theta + 4)\Phi + (4\theta^2 + 5\theta + 7)=0$.} 
& $(-\theta^2 - \theta - 1 , \pm (11\theta^2 + 13\theta + 17) )$ \\
\\
\hline
$0$ 
&
$(1, 3)-\infty$, 
&
$\infty$,\\
& 
$\displaystyle
\left(
\frac{\theta^2 + 2\theta + 1}{3}, 
\frac{10\theta^2 + 8\theta + 13}{3}
\right)-\infty
$ 
& 
$\displaystyle
\left(
\frac{\theta^2 + 2\theta + 1}{3} , \frac{\pm (10\theta^2 + 8\theta + 13)}{3} 
\right)
$, \\
& & $(1 ,\pm 3)$ \\ 
\\
\hline
$1$
&
$D_1=(-\theta^2 - \theta - 1, -40\theta^2 - 53\theta - 67) -\infty$,
&
$\infty$,\\ 
&
$D_2=(- 1, 3\theta + 3)-\infty$,
& 
$(-\theta^2 - \theta - 1 , \pm (40\theta^2 + 53\theta + 67))$, \\
&
$
\displaystyle
D_3=\sum_{i=5,6}(\Phi_i,(2\theta - 2)\Phi_i - \theta + 1)
-2\infty
$~\footnote{
$\Phi_5$, $\Phi_6$ are the roots of
$3\Phi^2 + (\theta^2 - \theta - 2)\Phi + (-2\theta^2 + 2\theta + 1)=0$. 
}, 
& $(-1 , \pm (3\theta + 3) )$ \\
\hline
$2$
&
$\emptyset$
&
$\infty$\\
\hline
\end{tabular}
\label{table:2310}
\end{minipage}
\end{table}

Note that in each case the rank $r \leq 3=d(g-1)$ where $d=[K:\Q]=3$
and $g=2$ is the genus. 

We implemented our method in {\tt MAGMA}. Our program 
succeeded in determining $C_s(K)$  for all $-2 \leq s \leq 2$,
and the results are given in Table~\ref{table:2310}. 
The entire computation
took approximately $2.5$ hours on a $2.8$ GHz Dual-Core AMD Opteron;
this includes the time taken for  
computing Selmer groups and  searching for points on the Kummer surfaces.
It is appropriate to give more details and we do this for the case $s=1$.
Let $C=C_1$ and write $J$ for its Jacobian. Let
\[
\cK=
\left\{ 
\infty,\; P_0,\; P_0^\prime,\; P_1,\; P_1^\prime
\right\},
\]
where
\[
P_0=(-\theta^2 - \theta - 1 , 40\theta^2 + 53\theta + 67), \qquad
P_1= (-1 , 3\theta + 3),  
\]
and $P_0^\prime$, $P_1^\prime$ are respectively the images of $P_0$,
$P_1$ under the hyperelliptic involution.
Let $D_1$, $D_2$, $D_3$ be the basis given in Table~\ref{table:2310}
for a subgroup of $J(K)$ of finite index. 
Let $L_0=\langle D_1, D_2, D_3 \rangle$. Our program verified that the
index of $L_0$ in $J(K)$ is not divisible by any prime $<75$. 
Our program used the point 
\[
P_0=(-\theta^2 - \theta - 1 , 40\theta^2 + 53\theta + 67)
\]
as the base point for the Abel-Jacobi map $\jmath$. The image of $\cK$ under $\jmath$ is
\[
\jmath(\cK)=\{D_1, \quad 0, \quad  2D_1, \quad D_1 + D_2, \quad D_1 - D_2\},
\]
where we have listed the elements of $\jmath(\cK)$ so that they correspond
to the above list of points of $\cK$.
Next our program applied the Mordell--Weil sieve as in Lemma~\ref{lem:mw}.
The program chose $22$ places $\vv$ which are places of good reduction for $C$
and such that $\#J(k_\vv)$ is divisible only by primes $<75$. In the notation
of Lemma~\ref{lem:mw},
\[
L_{22}=\langle 1386000 D_1 + 16632000 D_2 + 18018000 D_3, \quad
    24948000 D_2, \quad
    24948000 D_3 \rangle , 
\]
and
\begin{gather*}
W_{22}=\{
    0, \;
    D_1 - D_2, \;
    D_1, \;
    D_1 + D_2, \;
    2 D_1, \;
    D_1 + 12474000 D_2 + 87318000 D_3, \\
    277201 D_1 + 5821200 D_2 + 51004800 D_3,  \quad
    277201 D_1 - 6652800 D_2 - 36313200 D_3, \\
    -277199 D_1 + 6652800 D_2 + 36313200 D_3, \quad
    -277199 D_1 - 5821200 D_2 - 51004800 D_3
\}.
\end{gather*}

Next we would like to apply Theorem~\ref{thm:mwc} and so we need primes $p$
satisfying conditions (a)--(d) of that theorem. In particular, our program
searches for odd primes $p$, unramified in $K$, so that every place $\vv \mid p$
is a place of good reduction for $C$, and $\#J(k_\vv)$ is divisible only by
primes $<75$,
and so that $L_{22}$ is contained in the 
kernel of the homomorphism \eqref{eqn:proker}.
The smallest prime satisfying these conditions is $p=109$ which splits
completely in $K$ and so there are three degree $1$ places $\vv_1$,
$\vv_2$, $\vv_3$ above $109$. It turns out that
\[
J(k_\vv) \cong (\Z/110)^2,
\]
for $\vv=\vv_1,\vv_2,\vv_3$. The reader can easily see that
\[
L_{22} \subset 110 L_0 \subseteq 110 J(K)
\]
and so clearly $L_{22}$ is in the kernel of \eqref{eqn:proker} with $p=109$.
Moreover, the reader will easily see that every $\mathbf{w} \in W_{22}$
is equivalent modulo $110 L_0$ to some element of $\jmath(\cK)$. Hence
conditions (a), (c), (d) of the Theorem~\ref{thm:mwc}
are satisfied for each $\mathbf{w} \in W_{22}$ with $p=109$.
To show that $C(K)=\cK$ it is enough to show that $\tilde{M}_{109}(Q)$
has rank $3$ for all $Q \in \cK$.

It is convenient to take
\[
\omega_1=\frac{dx}{y}, \qquad \omega_2=\frac{xdx}{y},
\]
as basis for the $1$-forms on $C$. With this choice we computed the matrices
$\tilde{M}_{109}(Q)$ for $Q \in \cK$. For example, we obtained
\[
\tilde{M}_{109}(\infty)=
\begin{pmatrix}
 79 & 64 &  0\\
 31 &  0 &  0\\
104 &  0 & 82\\
\end{pmatrix} \pmod{109};
\]
this matrix of course depends on our choice of 
$U$ used to compute the Hermite Normal Form on page \pageref{HNF},
though as observed in the remarks after Theorem~\ref{thm:1},
its rank is independent of this choice of $U$.
The matrix $\tilde{M}_{109}(\infty)$ 
clearly has non-zero determinant and so rank $3$. It turns out that 
the four other $\tilde{M}_{109}(Q)$ also have rank $3$. This completes
the proof that $C(K)=\cK$.

We now return to the general case where $-2 \leq s \leq 2$, and would like to
recover the coprime integer solutions $u$, $v$ to equation~\eqref{eqn:335}
from the $K$-rational points on
$C_s$
and hence the solutions $(x,y,z)$
to \eqref{eqn:2310} with $y$ even and $x \equiv z^5 \pmod{4}$. 
From \eqref{eqn:XY} and \eqref{eqn:uvalpha} we see that
\[
Y=\frac{3(u+v\theta)}{\alpha^5}=3 \epsilon^s 
\left(\frac{u+v\theta}{u-v\theta} \right).
\]
Thus
\[
\frac{u}{v}=\theta \cdot \left( \frac{Y+3\epsilon^s}{Y-3\epsilon^s}\right).
\]
Substituting in here the values of $Y$ and $s$ from the $K$-rational points on the curves $C_s$, the 
only {\bf $\Q$-rational} values for $u/v$ we obtain are
respectively
$-1$, $2$, $0$, $5/4$, $1$; these respectively come from the points
$(\theta^2 + \theta + 1 , -\theta^2 - 2\theta - 1)$, 
$ (-\theta^2 - \theta - 1 , -11 \theta^2 - 13\theta - 17 )$,
$(1 , -3 )$,
$(-\theta^2 - \theta - 1 , 40 \theta^2 + 53\theta + 67 )$,
$(-1 , 3\theta + 3 )$. This immediately allows us to complete the proof of Theorem~\ref{thm:Ferm}.

The reader can find the {\tt MAGMA} code for verifying the
above computations at:   
{\tt http://www.warwick.ac.uk/staff/S.Siksek/progs/chabnf/}

\bigskip

\noindent {\bf Remarks.} 
\begin{enumerate}
\item[(i)]
Although our approach solves equation~\eqref{eqn:2310} completely, 
we point out that it is possible to eliminate some cases by using Galois representations
and level-lowering as Dahmen \cite{Dahmen} does for the equation $x^2+z^{10}=y^3$.
Indeed, by mimicking Dahmen's approach and making use of the work of Darmon and Merel \cite{DM},
and the so called \lq method for predicting the exponents of constants\rq\ \cite[Section 15.7]{Cohen} we were able to reduce
to the case $s=1$, and it is this case that corresponds to our non-trivial
solution $(x,y,z)=(\pm 3,-2,\pm 1)$. It seems however that the approach via Galois representations 
cannot in the current state of knowledge deal with case $s=1$.
\item[(ii)] Note that to solve our original problem \eqref{eqn:2310},
we did not need all $K$-rational points on the curves $C_s$,
merely those $(X,Y) \in C_s(K)$ with 
\[
\theta \cdot \left(\frac{Y+3 \epsilon^s}{Y-3 \epsilon^s}\right) \in \Q.
\]
This suggests that a higher dimensional analogue of 
elliptic curve Chabauty \cite{Br1}, \cite{Br2},
\cite{FW1}, \cite{FW2} should be applicable,
and this should give an alternative approach to \eqref{eqn:2310}.
Although it was not needed here we expect that this
idea will be useful in other contexts. 
\end{enumerate}

\end{document}